\documentclass[12pt,reqno]{amsart}
\usepackage{amsmath,amsthm}
\usepackage{amssymb}
\usepackage[english]{babel}
\usepackage{times}
\usepackage[margin=1.5in]{geometry}
\usepackage{mathrsfs}
\usepackage{comment}

\usepackage{graphicx}
\usepackage{tikz}
\usetikzlibrary {automata,positioning}

\newcommand{\z}{{\mathbb Z}}

\newcommand{{\rr}}{\mathcal{R}}

\newtheorem{thm}{Theorem}[section]
\newtheorem{corr}[thm]{Corollary}

\newtheorem{prop}[thm]{Proposition}

\newtheorem{que}[thm]{Question}

\newcommand{\pf}{\noindent\textbf{Proof.}$\ $}
\newcommand{\zb}{$\hfill\Box$}

\makeatletter\@addtoreset{equation}{section} \makeatother

\begin{document}
%\mag=1200

\title[Spectral dynamics for the Basilica group]{Spectral dynamics for the Basilica group}

\author[H. Liu]{Hongyi Liu}
\footnotemark[1]
\address{Hongyi Liu: School of Mathematics, Southeast University,
Nanjing, Jiangsu 211189, China.} \email{hongyiliu@seu.edu.cn}

\author[W. He]{Wei He}
\footnotemark[1]
\address{Wei He (corresponding author): School of Mathematics, Southeast University,
Nanjing, Jiangsu 211189, China.} \email{hewei@seu.edu.cn}

%\maketitle

\begin{abstract}
The projective  spectrum of a tuple $A=(A_0,A_1,\dots,A_n)$ with elements in a unital Banach algebra $\mathcal{A}$ is the collection of $[z_0:z_1:\cdots:z_n]\in \mathbb{P}^n$ such that $z_0A_0+z_1A_1+\dots+z_nA_n$ is not invertible in $\mathcal{A}$. For the tuple $A_{\rho}=(\rho({\bf e}),\rho({\bf a}),\rho({\bf b}))$, where ${\bf e},{\bf a},{\bf b}$ are the three states of the automaton generating the Basilica group $\mathcal{B}$ and $\rho$ is the Koopman representation, a dynamical map $F$ preserving the projective spectrum of $A_{\rho}$ is established. Further, a substantial common subset of the Julia set of $F$ and the projective spectrum of $A_{\rho}$, as well as a large common subset of the Fatou set of $F$ and the projective resolvent set of $A_{\rho}$ is identified. The result provides partial support for the conjecture of equality of the Julia set of $F$ and the projective spectrum of $A_{\rho}$.
\\%标蓝是因为此处的空行


%47A13(1991–now)Several-variable operator theory (spectral, Fredholm, etc.)
%37F10(2000–now)Dynamics of complex polynomials, rational maps, entire and meromorphic functions; Fatou and Julia sets [See also 32A10, 32A20, 32H02, 32H04]
%20E08(1991–now)Groups acting on trees [See also 20F65]

\end{abstract}

\maketitle

\section{Introduction}

Let $\mathcal{A}$ be a unital Banach algebra and let $A=(A_0,A_1,\dots,A_n)$ be a tuple of elements in $\mathcal{A}$. For $z=(z_0,z_1,\dots,z_n)\in \mathbb{C}^{n+1}$, define
\[
A(z)=z_0A_0+z_1A_1+\cdots+z_nA_n.
\]
The \emph{projective spectrum} of $A$, introduced by Yang \cite{Ya1}, is defined as
\begin{equation*}
	P(A)=\left\{ z\in \mathbb{C}^{n+1} : A(z)\text{ is not invertible in }\mathcal{A} \right\}.
\end{equation*}
The \emph{projective resolvent set} is its complement $P^c(A)=\mathbb{C}^{n+1}\setminus P(A)$. Since $A(z)$ is homogeneous, both the projective spectrum and the projective resolvent set can naturally be viewed as subsets of the $n$-dimensional projective space $\mathbb{P}^n$. Let $\phi: \mathbb{C}^{n+1}\setminus\{0\}\to \mathbb{P}^n$ be the canonical quotient map. The sets
\[
p(A)=\phi(P(A)\setminus\{0\}), \qquad p^c(A)=\phi(P^c(A))
\]
are still called the projective spectrum and projective resolvent set of $A$, respectively. The projective spectrum provides a powerful tool for studying the joint action of the tuple and the interactions among its elements, and has therefore found applications in many fields.

A case of particular interest is when the tuple $A$ arises from a generating set of a finitely generated group. To be precise, let $G$ be a finitely generated group with a generating set $S=(g_1,\dots,g_s)$, and let $\pi$ be a unitary representation of $G$ on a complex Hilbert space $\mathcal{H}$. Consider
\[
A_{\pi}(z)=z_0I+z_1\pi(g_1)+\cdots+z_s\pi(g_s).
\]
Then the projective spectrum $P(A_{\pi})$ encodes information about both $G$ and the representation $\pi$. For instance, when $\pi$ is finite-dimensional, $P(A)$ coincides with the zero locus of the \emph{characteristic polynomial} $Q_{\pi}(z)=\det(A_{\pi}(z))$.
This topic dates back to Dedekind and Frobenius \cite{Cu,De,Fr} and has laid the foundation for the theory of group representations. In recent years, Hu and Yang \cite{HY} revealed a connection between the characteristic polynomial and irreducibility of the group  representation, while Liu and Wang \cite{LW} used characteristic polynomials to study the combinatorial structure of finite Coxeter groups. We refer the reader to the monograph \cite{Ya3} and the references therein for further details.

Another important direction concerns the spectral theory of self-similar groups. For example, Grigorchuk and his collaborators investigated the spectral theory of the Grigorchuk group $\mathcal{G}$ of intermediate growth \cite{GN,GNS1,GNS2,GS}. Among many other results, they discovered that $\mathcal{G}$ admits a self-similar representation on the rooted binary tree. More remarkably, this self-similarity induces a rational map on a certain spectral set of $\mathcal{G}$, whose dynamical properties are intimately linked to the spectral properties of the group. This interplay between dynamics and spectral theory is often referred to as \emph{spectral dynamics} \cite{GoY}. In this direction, Yang and his collaborators examined the case of the infinite dihedral group \cite{GoY,GrY,Ya2,Ya3}. They eventually proved that the projective spectrum of a generating set of the infinite dihedral group is exactly the Julia set of the corresponding dynamical map on the spectrum. This constitutes one of the rare instances in which the Julia set of a dynamical map can be completely computed. Further contributions along this line include \cite{ZYL}, in which the spectral dynamics of another tuple in the infinite dihedral group was studied. In the same paper, the authors also investigated the spectral dynamics of a particular tuple of the lamplighter group and proved that its Julia set is contained in the projective spectrum.

The infinite dihedral group and the lamplighter group are both generated by $(3,2)$-automata (meaning automata with 3 states and 2 alphabets). According to the classification by Grigorchuk and his collaborators \cite{BGK1,BGK2,BGK3}, $(3,2)$-automata can generate at most 122 non-isomorphic groups. However, not all of them admit spectral dynamics; even when they do, the resulting spectral dynamics may differ significantly from one another.

In this paper, among the self-similar groups generated by $(3,2)$-automata, we identify another self-similar group that admits spectral dynamics, namely the Basilica group $\mathcal{B}$. This group was introduced by Grigorchuk and Żuk in \cite{GZ1} and some of its spectral properties were studied in \cite{GZ2}. The Basilica group possesses a number of remarkable properties. For instance, it is a torsion-free, contracting, regular weakly branch group of exponential growth. Meanwhile, it is the first example of a group that is amenable but not subexponentially amenable \cite{BV}. The amenability of self-similar groups has since been studied further in a number of works \cite{AAV1,AAV2,BKN,Ka}. The Schreier graphs associated with the Basilica group have also attracted considerable attention, see for example \cite{BGJRT,CDD,DDMN,DDS,SB}.

In the present paper, we investigate the spectral dynamics of the Basilica group. Let $\{{\bf e},{\bf a},{\bf b}\}$ (where ${\bf e}$ means identity) be the three states in the automaton of the Basilica group $\mathcal{B}$, and let $\rho$ be the Koopman representation. We focus on the tuple $A_{\rho}=(\rho({\bf e}),\rho({\bf a}),\rho({\bf b}))$. By exploiting the self-similarity of $\mathcal{B}$, we obtain a map
\begin{equation}\label{eqmain}
	F([z_0:z_1:z_2])=\left[(z_0+z_2)z_0:(z_0+z_2)z_2:-z_1^2\right].
\end{equation}
Let $E(F)$ be the extended indeterminacy set of $F$. Then $F$ is a map from $\mathbb{P}^2\setminus E(F)$ to itself, preserving both the projective spectrum and the projective resolvent set of $A_{\rho}$.
For $n\ge 0$, denote the $n$-th iterate of $F$ as $F^{n}$, where $F^0$ is the identity map. Let $\mathcal{F}(F)$ be the Fatou set of $F$ and $\mathcal{J}(F)$ its Julia set. The following two theorems are the main results of this paper.

\begin{thm}\label{main_theorem1}
	Let  
	\emph{(i)} $R_0=\{[z_0:z_1:z_2]:|z_0|>|z_1|+|z_2|\}$, \\
	\emph{(ii)} $R_1=\{[z_0:z_1:z_2]:|z_1|>|z_0|+|z_2|\}$, \\ 
	\emph{(iii)} $R_2=\{[z_0:z_1:z_2]:|z_2|>|z_0|+|z_1|\}$.\\
	Then 
	\[
	\bigcup^{\infty}_{n=0} F^{-n}\left(\bigcup_{k=0}^2R_k\right)\subseteq \mathcal{F}(F)\cap p^c(A_{\rho}).
	\]
\end{thm}

\begin{thm}\label{main_theorem2}
	Let $L=\left\{[z_0:z_1:z_2]:z_0+z_1+z_2=0\right\}$. Then  
	\[
	\overline{\bigcup^{\infty}_{n=0} F^{-n}(L)}\subseteq \mathcal{J}(F)\cap p(A_{\rho}).
	\]
\end{thm}

In Section \ref{section_dynamic_F}, we show how the dynamical map $F$ in \eqref{eqmain} arises from the Basilica group and compute its extended indeterminacy set. In Section \ref{section_Fatou_resolvent} and Section \ref{section_Julia_spectrum}, we prove Theorem \ref{main_theorem1} and Theorem \ref{main_theorem2}, respectively.
Section 5 concludes the paper and outlines some open questions for future work.\\

\section{The spectral dynamics of the Basilica group}\label{section_dynamic_F}

In this section, we first explain how the dynamical map $F$ in \eqref{eqmain} arises from the Basilica group.

\begin{figure}[htbp]
\centering
\begin{tikzpicture}[shorten >=1pt,node distance=2cm,auto]
  \node[state, label=above left:${\bf e}$]  (e) {$id$};
  \node[state, label=right:${\bf b}$]  (b) [above right=of e] {$id$};
  \node[state, label=above right:${\bf a}$]  (a) [below right=of b] {$\sigma$};

  \path[->] (e) edge [loop below, looseness=8] node {0,1} (e)
            (a) edge [bend right=20]  node [swap]       {0} (b)
			(a) edge  node      {1} (e)
			(b) edge              node        {0} (e)
			(b) edge              node  [swap]      {1} (a);
\end{tikzpicture}
\caption{Automaton of the Basilica group}
\label{automaton_Basilica}
\end{figure}
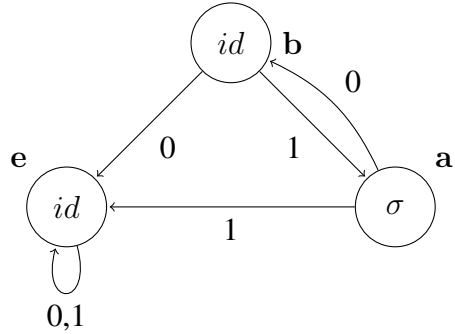

Recall that the Basilica group $\mathcal{B}$ is a self-similar group generated by the automaton in Figure \ref{automaton_Basilica}. It is a $(3, 2)$-automaton with $3$ states $\left\{{\bf e}, {\bf a},{\bf b}\right\}$ over a $2$-letter alphabet $\{0,1\}$. Then $\mathcal{B}$ admits a faithful action on the regular rooted binary tree $T$ in the following way. Let $0w$ and $1v$ be vertices of $T$. The action of the generators of $\mathcal{B}$ on $T$ is given by
\begin{equation*}\begin{aligned}
		&{\bf e}(0w)=0{\bf e}(w),&{\bf e}(1v)=1{\bf e}(v),\\
		&{\bf a}(0w)=1{\bf b}(w),&{\bf a}(1v)=0{\bf e}(v),\\	
		&{\bf b}(0w)=0{\bf e}(w),&{\bf b}(1v)=1{\bf a}(v).
	\end{aligned}
\end{equation*}
Let $\sigma$ be the permutation swapping $0$ and $1$,  
then the above action can be written simply as 
\begin{equation}\label{wreath_of_B}
	{\bf e}=({\bf e},{\bf e}),\, {\bf a}=\sigma ({\bf b},{\bf e}),\, 
	{\bf b}=({\bf e},{\bf a}).
\end{equation}

For a vertex $v$ of length $n$ in $T$, let $T_v$ be the subtree of $T$ with root $v$, and $\partial T_v$ be the boundary of $T_v$, i.e.,
\begin{equation*}
	\partial T_v=\left\{vv' : v'\text{ is an infinite sequence over } \left\{0,1\right\}\right\}.
\end{equation*}
The \emph{uniform Bernoulli measure} $\mu$ on $\partial T$ is defined by $\mu(\partial T_{v})=\frac{1}{2^n}$. Note that $T=T_{\varnothing}$, where $\varnothing$ is the empty word. Then $\mathcal{B}$ naturally admits a measure-preserving action on $\partial T$. 
Let $\mathcal{H} =L^2(\partial T,\mu)$. The \emph{Koopman representation} $\rho$ of $\mathcal{B}$ on $\mathcal{H}$ is defined as
\begin{equation*}
	(\rho (g)f)(w)=f(g^{-1}w),\qquad g\in \mathcal{B}, f\in \mathcal{H}, w\in \partial T.
\end{equation*}
Obviously, for any $g\in \mathcal{B}$, $\rho(g)$ is a unitary operator on $\mathcal{H}$, since the action of $g$ on $\partial T$ preserves the measure $\mu$. 

We consider the operator tuple $A_{\rho}=(\rho({\bf e}),\rho({\bf a}),\rho({\bf b}))$. Observe that $\rho({\bf e})=I$, and we simply denote the operators $\rho({\bf a})$ and $\rho({\bf b})$ as ${\bf a}$ and ${\bf b}$, respectively. Hence $A_{\rho}=(I,{\bf a},{\bf b})$. Then, (\ref{wreath_of_B}) implies that 
there exists a unitary operator $W:\mathcal{H}\rightarrow \mathcal{H}\oplus \mathcal{H}$ such that 
\begin{equation} \label{matrix_of_abc}	
	WIW^*=\left(\begin{matrix}
		I &O\\O &I
	\end{matrix}\right),
	W{\bf a}W^*=\left(\begin{matrix}
		O &I\\{\bf b} &O
	\end{matrix}\right),	
	W{\bf b}W^*=\left(\begin{matrix}
		I &O\\O &{\bf a}
	\end{matrix}\right). 
\end{equation}
Thus, we can define a function that preserves both the projective spectrum and the projective resolvent set of $A_{\rho}$.
This is the following proposition.

\begin{prop}\label{F_in_C_lemma}
	For $z=(z_0,z_1,z_2)\in \mathbb{C}^3$, $A_{\rho}(z)$ is invertible if and only if $A_{\rho}(\widetilde{F}(z))$ is invertible, where $\widetilde{F}:\mathbb{C}^3\rightarrow \mathbb{C}^3$ is defined as
	\begin{equation}\label{dynamic_in_C}
		\widetilde{F}(z)=\left((z_0+z_2)z_0,(z_0+z_2)z_2,-z_1^2\right).
	\end{equation}
\end{prop}
\pf By definition, $A_{\rho}(z)=z_0I +z_1{\bf a} +z_2{\bf b}$. In view of \eqref{matrix_of_abc}, we have 
\begin{equation}\label{Bz}
	A_{\rho}(z)=W^*\left(\begin{matrix}	(z_0+z_2)I & z_1 I\\z_1{\bf b} & z_0 I+z_2{\bf a}  \end{matrix}\right)W.
\end{equation}

We first consider the case when $z_1\neq 0$. By the argument of Schur complement (see \cite{Ga, GN, Sch}), we know that $A_{\rho}(z)$ is invertible if and only if 
\begin{equation*}
z_1{\bf b}-(z_0I+z_2{\bf a})(z_1I)^{-1}(z_0+z_2)I
\end{equation*}
is invertible. One verifies that the above formula is exactly  $-z_1^{-1}A_{\rho}(\widetilde{F}(z)).$
Thus, $A_{\rho}(z)$ is invertible if and only if $A_{\rho}(\widetilde{F}(z))$ is invertible. 

When $z_1=0$, it is easy to see from \eqref{Bz} that $A_{\rho}(z)$ is invertible if and only if both $(z_0+z_2)I$ and $z_0I+z_2{\bf a}$ are invertible. 
This holds if and only if $(z_0+z_2)z_0I+(z_0+z_2)z_2{\bf a}$ is invertible, which is exactly $A_{\rho}(\widetilde{F}(z))$ with $z_1=0$. The proof is complete.\zb\\

Observe that the map $\widetilde{F}:\mathbb{C}^3\rightarrow \mathbb{C}^3$ in Proposition \ref{F_in_C_lemma} preserves $P(A_{\rho})$ as well as $P^c(A_{\rho})$. To induce a map preserving $p(A_{\rho})$ and $p^c(A_{\rho})$, we first need to eliminate the preimage of the origin under $\widetilde{F}$. Recall that $\phi$ is the canonical map from $\mathbb{C}^3\setminus\{0\}$ to $\mathbb{P}^2$. For $S\subseteq \mathbb{C}^3$, we use the notation $\phi(S)$ to mean $\phi(S\setminus\left\{0\right\})$ for convenience. Let $F:\mathbb{P}^2\setminus \phi(\widetilde{F}^{-1}(0))\rightarrow \mathbb{P}^2$ be defined by $F([z])=\phi \widetilde{F} \phi ^{-1}$. It is easy to verify that $F$ is well-defined. The set 
$I_1(F)=\phi(\widetilde{F}^{-1}(0))$ is called the \emph{indeterminacy set} of $F$. 

We adopt the following conventions on $\mathbb{P}^2$. For $i=0,1,2$, let \(U_i=\{[z_0:z_1:z_2]\in \mathbb{P}^2: z_i\neq 0\}\) denote the standard affine chart. The topology on \(\mathbb{P}^2\) is induced by the Fubini-Study metric $d$, which is locally equivalent to the Euclidean metric on $\mathbb{C}^2$.

To proceed with iteration of $F$, let $I_n(F)=\phi(\widetilde{F}^{-n}(0))$ for any $n\ge 1$. It is obvious that $I_n(F) \subseteq I_{n+1}(F)$. Let $I_{\infty}(F)=\cup^{\infty}_{n=1}I_n(F)$ and $E(F)=\overline{I_{\infty}(F)}$. The set $E(F)$ is called the \emph{extended indeterminacy set} of $F$. Then the map
\begin{align*}
	F([z_0:z_1:z_2])=[(z_0+z_2)z_0:(z_0+z_2)z_2:-z_1^2],
\end{align*}
as well as $F^n$, is well-defined on $\mathbb{P}^2\setminus E(F)$. This is exactly the dynamical map in \eqref{eqmain}. 
Moreover, the set $E(F)$ can be explicitly computed.

\begin{prop}\label{indeterminacy_set}
	The extended indeterminacy set of $F$ in \eqref{eqmain} has the form 
	\begin{equation*}
		E(F)=\left\{[1:a:0] : \left\lvert a\right\rvert=1 \right\}\cup\left\{[1:0:a] : \left\lvert a\right\rvert=1 \right\}.
	\end{equation*}
\end{prop}

\pf Direct computation shows that $\widetilde{F}^{-1}(0)=\left\{(z,0,-z) : z\in \mathbb{C}\right\}$. Hence 
\begin{equation*}
	I_1(F)=\left\{[1:0:-1]\right\}.
\end{equation*}
For $n\geq 2$, since 
\begin{equation*}
	\widetilde{F}^{-n}(0)=\widetilde{F}^{-(n-1)}(0)\cup \widetilde{F}^{-1}\left(\widetilde{F}^{-(n-1)}(0)\setminus\widetilde{F}^{-(n-2)}(0)\right),
\end{equation*}
we have
\begin{align}\label{computing_Ik}
		I_n(F)=\phi(\widetilde{F}^{-n}(0))=I_{n-1}(F)\cup F^{-1}(I_{n-1}(F)\setminus I_{n-2}(F)),
\end{align}
where $I_0(F)$ stands for the empty set $\varnothing$. Observe that for $a\in \mathbb{C}$,
\begin{align}
	F^{-1}([1:0:a])&=\left\{[1:b:0]:b^2=-a\right\},\label{Ik_odd2even}\\
	F^{-1}([1:a:0])&=\begin{cases}
	    \left\{[1:0:a]\right\}, &\text{if } a\neq -1,\\
		\varnothing, &\text{if } a= -1.
		\end{cases}\label{Ik_even2odd}
\end{align}
Then we have the following two claims.
\par\textbf{Claim 1.} For $n\geq 1$, $I_{2n}(F)=I_{2n-1}(F)\cup \left\{[1:a:0]:a^{2^n}=1\right\}$. 

We prove Claim 1 by induction. For $n=1$, by (\ref{computing_Ik}) and (\ref{Ik_odd2even}), we know that
\begin{align*}
I_2(F)=I_1(F)\cup F^{-1}([1:0:-1])=I_1(F)\cup \left\{[1:a:0]:a^{2}=1\right\}.
\end{align*}
Assume that for $n=k$, $I_{2k}(F)=I_{2k-1}(F)\cup \left\{[1:a:0]:a^{2^k}=1\right\}$. Then (\ref{computing_Ik}) and (\ref{Ik_even2odd}) show that
\begin{align*}
	I_{2k+1}(F)&=I_{2k}(F)\cup F^{-1}\left\{[1:a:0]:a^{2^k}=1\right\}\\	
	&=I_{2k}(F)\cup \left\{[1:0:a]:a^{2^k}=1,\;a\neq-1\right\}.
	%&=I_{2k}(F)\cup \left\{[1:0:a]:a^{2^k}=1\right\}.
\end{align*}
But $I_1(F)=\left\{[1:0:-1]\right\}\subseteq I_{2k+1}(F)$, hence we have 
\begin{align*}%\label{Ikodd}
	I_{2k+1}(F)=I_{2k}(F)\cup \left\{[1:0:a]:a^{2^k}=1\right\}.
\end{align*}
Therefore, by (\ref{computing_Ik}) and (\ref{Ik_odd2even}),
\begin{equation*}
	\begin{aligned}
		I_{2(k+1)}(F)&=I_{2k+1}(F)\cup F^{-1}\left\{[1:0:a]:a^{2^k}=1\right\}\\
		&=I_{2k+1}(F)\cup \left\{[1:a:0]:a^{2^{(k+1)}}=1\right\},
	\end{aligned}
\end{equation*}
and Claim 1 is proved. 

\par\textbf{Claim 2.} For $n\geq 1$, $I_{2n+1}(F)=I_{2n}(F)\cup \{[1:0:a]:a^{2^n}=1\}$. 

The proof of Claim 2 is similar to that of Claim 1. %\eqref{Ikodd} from induction hypothesis in Claim 1. 

Claim 1 and Claim 2 indicate that 
\begin{equation*}
	I_{\infty}(F)=\left(\bigcup_{n=1}^{\infty}\left\{[1:0:a]:a^{2^n}=1\right\}\right)\cup \left(\bigcup_{n=1}^{\infty}\left\{[1:a:0]:a^{2^n}=1\right\}\right).
\end{equation*} Since $\bigcup_{n=1}^{\infty}\left\{a:a^{2^n}=1\right\}$ is dense in the unit circle, the equivalence of the Fubini-Study metric and the Euclidean metric on the affine chart $U_0$ yields that $E(F)=\overline{I_{\infty}(F)}$ has the desired form, completing the proof. \zb\\

\section{The projective resolvent set and the Fatou set}\label{section_Fatou_resolvent}
In this section, we identify a large common subset of the projective resolvent set of $A_{\rho}$ and the Fatou set of $F$. To begin with, we recall the concepts of normal families and Fatou sets.

Let $X$, $Y$ be metric spaces and let $\Omega\subseteq X$. A family of continuous functions $\{f_{\lambda}:\Omega\to Y:\lambda\in\Lambda\}$ is said to be \emph{normal} if every sequence in the family contains a subsequence that converges uniformly on every compact subset of $\Omega$. A point $p\in \mathbb{P}^2\setminus E(F)$ is called a \emph{Fatou point} of $F$ if there exists a neighborhood $U\subset \mathbb{P}^2\setminus E(F)$ of $p$ such that $\{F^n:n\ge 1\}$ is a normal family on $U$. The set of Fatou points of $F$ is called the \emph{Fatou set} of $F$, denoted by $\mathcal{F}(F)$, or simply $\mathcal{F}$ when no confusion arises. Note that $\mathcal{F}$ is open.

We shall show that the following sets constitute the common part of the Fatou set of $F$ and the projective resolvent set of $A_{\rho}$.

\begin{thm}\label{Fatou_resolvent}
	Let  
	\emph{(i)} $R_0=\{[z_0:z_1:z_2]:|z_0|>|z_1|+|z_2|\}$, \\
	\emph{(ii)} $R_1=\{[z_0:z_1:z_2]:|z_1|>|z_0|+|z_2|\}$, \\ 
	\emph{(iii)} $R_2=\{[z_0:z_1:z_2]:|z_2|>|z_0|+|z_1|\}$.\\
	Then 
	\[
	\bigcup^{\infty}_{n=0} F^{-n}\left(\bigcup_{k=0}^2R_k\right)\subseteq \mathcal{F}(F)\cap p^c(A_{\rho}).
	\]
\end{thm}

This theorem is divided into Proposition \ref{resolvent_Basilica}, \ref{R_0_Fatou}, \ref{R_1_Fatou} and \ref{R_2_Fatou}. We will prove them separately. 

\begin{prop}\label{resolvent_Basilica}
	The set $\bigcup^{\infty}_{n=0} F^{-n}\left(\bigcup_{k=0}^2R_k\right)$ is contained in $p^c(A_{\rho})$.
\end{prop}
\pf Since $F$ preserves $p^c(A_{\rho})$, it suffices to show that $\bigcup_{k=0}^2R_k\subseteq p^c(A_{\rho})$. 

%Denote the norm of operator $A$ as $\left\lVert A\right\rVert$. 

For $z=[z_0:z_1:z_2]\in R_0$, since ${\bf a}$, ${\bf b}$ are both unitary operators and $|z_1|+|z_2|<|z_0|$, we have $\lVert z_1 {\bf a}+z_2{\bf b}\rVert <|z_0|$. It follows that $z_0I+z_1{\bf a}+z_2{\bf b}$ is invertible. Hence $R_0\subseteq p^c(A_{\rho})$. 

For $z=[z_0:z_1:z_2]\in R_1$, we have that $|z_1|>|z_0|+|z_2|$. It is obvious that $| z_1|>0$, and $|{z_0}/{z_1}|+|{z_2}/{z_1}|<1$. Hence $z_0I+z_1{\bf a}+z_2{\bf b}$ is invertible if and only if 
\begin{align*}
	\frac{1}{z_1}{\bf a}^{-1}(z_0I+z_1{\bf a}+z_2{\bf b})=I+\frac{z_0}{z_1}{\bf a}^{-1}+\frac{z_2}{z_1}{\bf a}^{-1}{\bf b}
\end{align*}
is invertible. Since ${\bf a}^{-1}$ and ${\bf a}^{-1}{\bf b}$ are both unitary operators and $|{z_0}/{z_1}|+|{z_2}/{z_1}|<1$, it holds that $\left\lVert \frac{z_0}{z_1}{\bf a}^{-1}+\frac{z_2}{z_1}{\bf a}^{-1}{\bf b}\right\rVert <1$, indicating that $I+\frac{z_0}{z_1}{\bf a}^{-1}+\frac{z_2}{z_1}{\bf a}^{-1}{\bf b}$ is invertible. Thus $z_0I+z_1{\bf a}+z_2{\bf b}$ is invertible, and $R_1\subseteq p^c(A_{\rho})$. 

For $R_2$, the proof is similar to that of $R_1$.

In summary, we have $\bigcup_{k=0}^2R_k\subseteq p^c(A_{\rho})$, completing the proof.\zb\\

\begin{prop}\label{R_0_Fatou}
	For $R_0=\{[z_0:z_1:z_2]:|z_0|>|z_1|+|z_2|\}$, $\bigcup^{\infty}_{n=0} F^{-n}(R_0)\subseteq \mathcal{F}$.
\end{prop}
\pf By definition, the Fatou set is invariant under $F^{-1}$, so it suffices to show that $R_0\subseteq \mathcal{F}$. For simplicity, rewrite $R_0=\{[1:z_1:z_2]:1>|z_1|+|z_2|\}$. 

The map $F$ in \eqref{eqmain}, restricted to $R_0$, can be written as 
\begin{equation}\label{renormalized_F}
	F([1:z_1:z_2])=\left[1:z_2:\frac{-z_1^2}{1+z_2}\right].
\end{equation}
For $z=[1:z_1:z_2]\in R_0$, since 
\begin{align*}
	|z_2|+\left|\frac{-z_1^2}{1+z_2}\right|\leq |z_2|+\frac{|z_1|^2}{1-|z_2|}<|z_2|+|z_1|<1,	
\end{align*}
we have $F(z)\in R_0$. Thus $F$ can be iterated in $R_0$. 
Let 
\begin{equation}\label{recursive_relation}
	z_n=\frac{-z_{n-2}^2}{1+z_{n-1}},\qquad n\geq 3.
\end{equation}
Then for $n\ge 1$, it follows by induction that
\begin{equation*}
	F^n([1:z_1:z_2])=[1:z_{n+1}:z_{n+2}].
\end{equation*}

For fixed $z=[1:z_1:z_2]\in R_0$, let $\varepsilon=(1-|z_1|-|z_2|)/3$, $a_1=|z_1|+\varepsilon$, and $a_2=|z_2|+\varepsilon$. Hence $a_1+a_2=1-\varepsilon<1$. Set 
\begin{align*}
	U=\{[1:w_1:w_2]:|w_1-z_1|<\varepsilon, |w_2-z_2|<\varepsilon\}.
\end{align*}
It is obvious that $U\subset R_0$. We will show that $\{F^n\}$ converges uniformly on $U$, so that $z$ is a Fatou point of $F$.

Let $w=[1:w_1:w_2]\in U$. By the recursive relation $w_n=\frac{-w_{n-2}^2}{1+w_{n-1}}$, one verifies that for the odd terms, either all are zero (if $w_1=0$) or none are zero (if $w_1\neq 0$). The same holds for the even terms. Thus we may assume $w_1w_2\neq 0$. It follows that $w_n\neq 0$ for any $n$. Since $F^{2n-2}(w)\in R_0$, we have 
$|w_{2n-1}|+|w_{2n}|<1$. Then the recursive relation \eqref{recursive_relation} implies that 
\begin{align*}
	\left|\frac{w_{2n+1}}{w_{2n-1}}\right|=\left|\frac{w_{2n-1}}{1+w_{2n}} \right|
	\leq \frac{|w_{2n-1}|}{1-|w_{2n}|}<1, 
\end{align*}
which means that $\{|w_{2n-1}|\}$ is monotonically decreasing. Same arguments imply that $\{|w_{2n}|\}$ is also monotonically decreasing. Then, since $|w_1|<|z_1|+\varepsilon=a_1$ and $|w_2|<|z_2|+\varepsilon=a_2$, we have
\begin{align*}
	\left|\frac{w_{2n+1}}{w_{2n-1}}\right|\leq \frac{\left|w_{2n-1}\right|}{1-|w_{2n}|}<\frac{|w_1|}{1-|w_2|}<\frac{a_1}{1-a_2}. 
\end{align*}
It follows that 
\begin{align*}
	|w_{2n+1}|<\left(\frac{a_1}{1-a_2}\right)^{n}|w_1|<\left(\frac{a_1}{1-a_2}\right)^{n}a_1.
\end{align*}
Consequently, $\{w_{2n-1}\}$ converges to $0$ uniformly on $U$. The same method also gives that $\{w_{2n}\}$ converges to $0$ uniformly on $U$. Since both odd and even subsequences of $\{w_n\}$ converge uniformly to $0$, we have $F^n(w)$ converges uniformly to the constant function $[1:0:0]$ on $U$ with respect to the Fubini-Study metric. Indeed, the Fubini-Study metric is equivalent to the Euclidean metric on $U_0$, so the uniform convergence of coordinates implies uniform convergence in the Fubini-Study metric. This implies that $z$ is a Fatou point, and the proof is complete. \zb\\

\begin{prop}\label{R_1_Fatou}
	For $R_1=\{[z_0:z_1:z_2]:|z_1|>|z_0|+|z_2|\}$, $\bigcup^{\infty}_{n=0} F^{-n}(R_1)\subseteq \mathcal{F}$.
\end{prop}
\pf Since the Fatou set is invariant under $F^{-1}$, it suffices to show that $R_1\subseteq \mathcal{F}$. Since $z_1\neq 0$ for $z\in R_1$, we can rewrite $R_1=\{[z_0:1:z_2]:1>|z_0|+|z_2|\}$. 

For $z=[z_0:1:z_2]\in R_1$, let $\varepsilon=(1-|z_0|-|z_2|)/3$, $a_0=|z_0|+\varepsilon$, and $a_2=|z_2|+\varepsilon$. It follows that $a_0+a_2=1-\varepsilon<1$. Set
\begin{equation*}
	U=\{[w_0:1:w_2]:|w_0-z_0|<\varepsilon, |w_2-z_2|<\varepsilon\}.
\end{equation*}
We show that both $\{F^{2n-1}\}$ and $\{F^{2n}\}$ converge uniformly on $U$.

On $\mathbb{P}^2\setminus E(F)$, the map $F$ acting on $w=[w_0:w_1:w_2]$ can be written in the piecewise form
\begin{equation}\label{F_R1}
	F(w)=
	\begin{cases}
		\left[w_0:w_2:\frac{-w_1^2}{w_0+w_2}\right], & \text{if } w_2\neq -w_0,\\[1.2ex]
		[0:0:1], & \text{if } w_2=-w_0.
	\end{cases}
\end{equation}	
Observe that on the chart $\{w_2\neq -w_0\}$, the first coordinate remains fixed, which greatly simplifies the iteration of $F$. When $w_2=-w_0$, one checks that $F^{2n-1}(w)=[0:0:1]$ and $F^{2n}(w)=[0:1:0]$ for $n\ge 1$. Thus the orbit of $w$ falls into the $2$-cycle $[0:0:1]\leftrightarrow[0:1:0]$. 

Now for $w=[w_0:1:w_2]\in U$, we first consider the case that the orbit $\{F^n(w)\}$ never hits the set $\{w_2=-w_0\}$. Let $w_1=1$ and 
\begin{align}\label{recursive_R1}
	w_n=\frac{-w_{n-2}^2}{w_0+w_{n-1}},\qquad n\ge 3.
\end{align}
By induction, we have for $n\ge 1$,
\begin{equation}\label{R_1_iteration}
	F^n([w_0:1:w_2])=[w_0:w_{n+1}:w_{n+2}].
\end{equation}
Since $|w_0|<|z_0|+\varepsilon=a_0$ and $|w_2|<|z_2|+\varepsilon=a_2$, applying \eqref{recursive_R1} gives  
\begin{equation*}
	|w_{3}|=\left|\frac{-w_1^2}{w_0+w_2}\right|\ge\frac{1}{|w_0|+|w_2|}>\frac{1}{a_0+a_2}>1,
\end{equation*}
and hence
\begin{equation*}
	|w_{4}|=\left|\frac{-w_2^2}{w_0+w_3}\right|\le\frac{|w_2|^2}{|w_3|-|w_0|}<\frac{a_2^2}{1-a_0}<a_2.
\end{equation*}
Repeating the above argument inductively, we have $|w_{2n+1}|> 1$ and $|w_{2n}|<a_2$, for all $n\ge 1$.
Thus, for the odd terms,
\begin{equation*}
	\left| \frac{w_{2n+1}}{w_{2n-1}}\right|
	=\frac{|w_{2n-1}|}{| w_0+w_{2n}|}
	\ge \frac{|w_{2n-1}|}{|w_0|+|w_{2n}|}
	>\frac{1}{a_0+a_2}>1.
\end{equation*}
Consequently, 
\begin{equation}\label{R_1_increasing_seq}
	|w_{2n+1}|>\left(\frac{1}{a_0+a_2}\right)^n|w_1|=\left(\frac{1}{a_0+a_2}\right)^n.
\end{equation}
For the even terms,
\begin{equation*}
	\left| \frac{w_{2n+2}}{w_{2n}}\right|
	=\frac{|w_{2n}|}{|w_0+w_{2n+1}|}
	\leq \frac{|w_{2n}|}{|w_{2n+1}|-|w_0|}
	<\frac{a_2}{1-a_0}<1.
\end{equation*}
It follows that 
\begin{equation}\label{R_1_decreasing_seq}
	|w_{2n+2}|<\left(\frac{a_2}{1-a_0}\right)^n\left| w_2\right|<\left(\frac{a_2}{1-a_0}\right)^n.
\end{equation}
Inserting these estimates into \eqref{R_1_iteration} and using the Fubini-Study metric, we see that $\{F^{2n-1}(w)\}$ converges uniformly to the constant map $[0:0:1]$ and $\{F^{2n}(w)\}$ converges uniformly to $[0:1:0]$ for those points in $U$ whose orbits avoid the set $\{w_2=-w_0\}$.

It remains to treat those points in $U$ whose orbits hit $\{w_2=-w_0\}$ at some finite time. For such a point $w\in U$, suppose $N(w)$ is the first positive integer such that $w_{N(w)+2}=-w_0$ (using the notation $F^n(w)=[w_0:w_{n+1}:w_{n+2}]$). When $n\le N(w)$, the estimates in \eqref{R_1_increasing_seq} and \eqref{R_1_decreasing_seq} imply that the convergence speed of the odd and even subsequences of $\{F^n(w)\}$ is independent of $w$. For $n> N(w)$, the orbit has entered the $2$-cycle $[0:0:1]\leftrightarrow[0:1:0]$ and thus the odd and even subsequences of $\{F^n(w)\}$ have reached the corresponding limit, respectively. 

The above discussions prove uniform convergence of the odd and even subsequences of $\{F^n(w)\}$ on the whole $U$. Consequently, $\{F^n\}$ is a normal family on $U$, so $z$ is a Fatou point. Since $z\in R_1$ was arbitrary, we have $R_1\subset \mathcal{F}$. The proof is complete.\zb\\

\begin{prop}\label{R_2_Fatou}
	For $R_2=\{[z_0:z_1:z_2]:|z_2|>|z_0|+|z_1|\}$, $\bigcup^{\infty}_{n=0} F^{-n}(R_2)\subseteq \mathcal{F}$.
\end{prop}
\pf For $z=[z_0:z_1:z_2]\in R_2$, it is obvious that $z_0+z_2\neq 0$. Hence $F(z)$ can be written as 
\begin{equation*}
	F(z)=\left[z_0:z_2:\frac{-z_1^2}{z_0+z_2}\right].
\end{equation*}
It follows that
\begin{align*}
	|z_0|+\left|\frac{-z_1^2}{z_0+z_2}\right|\leq |z_0|+\frac{|z_1|^2}{|z_2|-|z_0|}<|z_0|+|z_1|<|z_2|.	
\end{align*}
This implies that $F(z)\in R_1$. Consequently, Proposition \ref{R_1_Fatou} and the invariance of $\mathcal{F}$ under $F^{-1}$ yield the result.\zb\\

\section{The projective spectrum and the Julia set}\label{section_Julia_spectrum}
In this section, we identify a substantial common subset of the projective spectrum of $A_{\rho}$ and the Julia set of $F$.

Recall that the \emph{Julia set} of $F$, denoted by $\mathcal{J}(F)$, or simply $\mathcal{J}$, is defined as the complement of the Fatou set $\mathcal{F}(F)$. Points belonging to $\mathcal{J}(F)$ are called \emph{Julia points} of $F$. Note that $\mathcal{J}(F)$ is closed and contains the extended indeterminacy set $E(F)$. In this section, we characterize a substantial portion of the projective spectrum $p(A_{\rho})$ and demonstrate that this portion is also contained in $\mathcal{J}(F)$. Specifically, we prove the following main theorem.

\begin{thm}\label{Julia_spectrum}
	Let $L=\left\{[z_0:z_1:z_2]:z_0+z_1+z_2=0\right\}$. Then  
	\[
	\overline{\bigcup^{\infty}_{n=0} F^{-n}(L)}\subseteq \mathcal{J}(F)\cap p(A_{\rho}).
	\]
\end{thm}

The proof of Theorem \ref{Julia_spectrum} is divided into two parts. The first part (Proposition \ref{L_in_spectrum}) proves that $\overline{\bigcup^{\infty}_{n=0} F^{-n}(L)}\subseteq  p(A_{\rho})$, and the second part (Proposition \ref{L_in_Julia}) proves that $\overline{\bigcup^{\infty}_{n=0} F^{-n}(L)}\subseteq \mathcal{J}$.

To prove the first part, we use the method of operator recursions. Recall that for the binary tree $T$, we have defined $\mathcal{H} =L^2(\partial T,\mu)$, where $\mu$ is the uniform Bernoulli measure on $\partial T$. For $n\ge 0$, let $\mathcal{H}_{n}$ be the subspace of $\mathcal{H}$ spanned by $\left\{\chi_w : w\in \left\{0,1\right\}^n\right\}$, where $\chi_w$ is the characteristic function of $\partial T_w$. Then $\mathrm{dim}\mathcal{H}_{n}=2^n$.  Since every $\chi_w$, where $w\in \left\{0,1\right\}^n$ and $n\ge 1$, is of the form $\chi_{0v}$ or $\chi_{1v}$ for some $v\in \{0,1\}^{(n-1)}$, we obtain a natural isomorphism between $\mathcal{H}_{n}$ and $\mathcal{H}_{n-1}\oplus \mathcal{H}_{n-1}$. 

Note that for each $n\ge 0$, $\mathcal{H}_{n}$ is invariant under $\rho(g)$ for any $g\in\mathcal{B}$. Let $\rho_n=\rho|_{\mathcal{H}_{n}}$. Set $I_n=\rho_n({\bf e})$, ${\bf a}_n=\rho_n({\bf a})$ and ${\bf b}_n=\rho_n({\bf b})$. Since $\mathcal{H}_{n}$ is finite dimensional, $I_n$, ${\bf a}_n$ and ${\bf b}_n$ are all matrices. By appropriately arranging the order of the basis of  $\mathcal{H}_{n}$, it holds that
\begin{equation}
	I_n=\left(\begin{matrix}
		I_{n-1}&O\\O&I_{n-1}
	\end{matrix}\right),\, 	{\bf a}_n=\left(\begin{matrix}
		O&I_{n-1}\\{\bf b}_{n-1} &O
	\end{matrix}\right),\,	{\bf b}_n=\left(\begin{matrix}
		I_{n-1}&O\\O&{\bf a}_{n-1}
	\end{matrix}\right),
	\label{matrix_of_abc_finite}
\end{equation}
for $n\geq 1$, where $O$ is the zero matrix and $I_0={\bf a}_0={\bf b}_0=1$. Thus, for any $n$, $I_n$ is the identity matrix and ${\bf a}_n$ and ${\bf b}_n$ are unitary matrices. 

For $n\ge 0$, let $A_{\rho_n}=(I_n, {\bf a}_n, {\bf b}_n)$. Denote $A_{\rho_n}(z)=z_0I_n+z_1{\bf a}_n+z_2{\bf b}_n$ for $z=(z_0,z_1,z_2)\in \mathbb{C}^3$. For a fixed $z$, it is clear that $A_{\rho}(z)$ is not invertible in $B(\mathcal{H})$ if there exists $n\ge 0$ such that $A_{\rho_n}(z)$ is not invertible in $B(\mathcal{H}_{n})$. Therefore, 
\begin{equation*}
	P(A_{\rho})\supseteq \bigcup_{n=0}^{\infty} P(A_{\rho_n}).
\end{equation*}
Now we are ready to prove the following proposition.
\begin{prop}\label{L_in_spectrum}
	The set $\overline{\bigcup^{\infty}_{n=0} F^{-n}(L)}$ is contained in $p(A_{\rho})$. 
\end{prop}
\pf In view of  (\ref{matrix_of_abc_finite}), we have 
\begin{equation*}
	\mathrm{det}(A_{\rho_n}(z))=\mathrm{det}\left(\begin{matrix}
		(z_0+z_2)I_{n-1} &z_1I_{n-1}\\z_1{\bf b}_{n-1} &z_0I_{n-1}+z_2{\bf a}_{n-1}
	\end{matrix}\right).
\end{equation*}
Using the Schur complement arguments as in the proof of Proposition \ref{F_in_C_lemma}, we know that $\mathrm{det}(A_{\rho_n}(z))=\mathrm{det}(A_{\rho_{n-1}}(\widetilde{F}(z)))$. Using this relation repeatedly, we have that $\mathrm{det}(A_{\rho_n}(z))=\mathrm{det}(A_{\rho_0}(\widetilde{F}^n(z)))$. Thus $A_{\rho_n}(z)$ is not invertible, if and only if $\mathrm{det}(A_{\rho_n}(z))=0$, if and only if $\mathrm{det}(A_{\rho_0}(\widetilde{F}^n(z)))=0$. It is easy to compute that  $\mathrm{det}(A_{\rho_0}(z))=z_0+z_1+z_2$. Hence if we let $L'=\left\{(z_0,z_1,z_2):z_0+z_1+z_2=0\right\}$, then $A_{\rho_n}(z)$ is not invertible if and only if $\widetilde{F}^n(z)\in L'$. Thus, $P(A_{\rho_n})=\widetilde{F}^{-n}(L')$. Consequently, 
\begin{equation*}
	P(A_{\rho})\supseteq \bigcup_{n=0}^{\infty} \widetilde{F}^{-n}(L').
\end{equation*}
Let $\phi$ be the canonical map from $\mathbb{C}^3\setminus \{0\}$ to $\mathbb{P}^2$. Then
\begin{equation*}
	p(A_{\rho})\supseteq \bigcup^{\infty}_{n=0} F^{-n}\left(\phi(L')\right)=\bigcup^{\infty}_{n=0} F^{-n}\left(L\right).
\end{equation*}
Since $p(A_{\rho})$ is closed \cite{Ya1}, we finally obtain
\begin{equation*}
p(A_{\rho})\supseteq \overline{\bigcup^{\infty}_{n=0} F^{-n}\left(L\right)}.
\end{equation*}
which is the desired result. \zb\\

The proof of Proposition \ref{indeterminacy_set} indicates that the extended indeterminacy set satisfies $E(F)=\overline{\bigcup^{\infty}_{n=0}F^{-n}(I_1)}$, where $I_1=\left\{[1:0:-1]\right\}$. Since $I_1\subseteq L$, we have  $E(F)\subseteq \overline{\bigcup^{\infty}_{n=0} F^{-n}(L)}\subseteq p(A_{\rho})$. This yields the following corollary.
\begin{corr}\label{corr_E_in_p}
	The extended indeterminacy set $E(F)$ is contained in $p(A_{\rho})$.
\end{corr}
The spectra of ${\bf a}$ and ${\bf b}$ can be derived from Corollary \ref{corr_E_in_p}. 
\begin{corr}
	The spectra of ${\bf a}$ and ${\bf b}$ are both the unit circle. 
\end{corr}
\pf
Since $\left\{[1:a:0]: |a|=1 \right\}\subseteq E(F)\subseteq p(A_{\rho})$, we have $\left\{(1,a,0): |a|=1 \right\}$ is contained in $P(A_{\rho})$. Hence the definition of the projective spectrum shows that
the spectrum of ${\bf a}$ contains the unit circle. 
On the other hand, ${\bf a}$ is unitary, hence its spectrum is contained in the unit circle. Therefore, the spectrum of ${\bf a}$
is exactly the unit circle. The same argument applies verbatim to ${\bf b}$. This completes the proof. \zb\\

We now turn to the proof of the second part of Theorem \ref{Julia_spectrum}, which constitutes the most delicate part of  the proof. 

\begin{prop}\label{L_in_Julia}
	The set $\overline{\bigcup^{\infty}_{n=0} F^{-n}(L)}$ is contained in $\mathcal{J}$. 
\end{prop}

\pf To show that $\overline{\bigcup^{\infty}_{n=0} F^{-n}(L)}\subseteq \mathcal{J}$, it suffices to prove $L\subseteq \mathcal{J}$, since $\mathcal{J}$ is closed and invariant under $F^{-1}$. We will prove the following two cases.

\textbf{Case 1.} If $z$ is a point in $L$ with all its coordinates nonzero, then $z\in\mathcal{J}$.

Recall that $L=\{[z_0:z_1:z_2]:z_0+z_1+z_2=0\}$. Any point in $L$ with all coordinates nonzero can be written uniquely as
$z_a=[1:a:-a-1]$ for some $a\in\mathbb{C}\setminus\{0,-1\}$. We prove $z_a\in\mathcal{J}$.

Assume, for contradiction, that $z_a$ is a Fatou point of $F$. A direct computation gives $F(z_a)=[1:-a-1:a]$ and $F^2(z_a)=z_a$. Set $G=F^2$. Then $z_a$ is a fixed point of $G$. Since $\{G^n\}$ is a subsequence of $\{F^n\}$, $z_a$ is also a Fatou point of $G$. Hence  there exists $r>0$ such that $B(z_a,r)=\{d(z,z_a)<r\}\subset U_0\setminus E(F)$, and the family $\{G^n\}$ is normal on $B(z_a,r)$. By Arzel\`a-Ascoli's theorem, $\{G^n\}$ is equicontinuous on every compact subset of $B(z_a,r)$.

For this $r$, equicontinuity of $\{G^n\}$ gives some $t\in(0,r)$ such that for all $z$ with $d(z,z_a)\le t$ and all $n>0$,
\[
d(G^n(z),G^n(z_a))=d(G^n(z),z_a)<r.
\]
Let $K=\overline{B(z_a,t)}$. Then $G^n(K)\subseteq B(z_a,r)\subset U_0\setminus E(F)$ for all $n$. Moreover, since $K$ is a compact subset of $B(z_a,r)$, $\{G^n\}$ is equicontinuous on $K$.

Now we transfer the dynamics to local coordinates centered at $z_a$. For $w=(w_1,w_2)^{\mathsf{T}}\in\mathbb{C}^2$, define $\varphi(w)=[1:w_1+a:w_2-a-1]$. Set
\[
K_0=\bigl\{w\in\mathbb{C}^2:
\varphi(w)\in K\bigr\}.
\]
Since $\varphi$ is a homeomorphism from $K_0$ to $K$, $K_0$ is compact and contains $\mathbf{0}$ as an interior point. Define $\hat{G}$ on $K_0$ by
\[
\hat{G}=\varphi^{-1}G\varphi.
\]
Since $G=F^2$, a direct computation yields		
\begin{align*}
	\hat{G}(w)
	=\begin{pmatrix}
		-\dfrac{(w_1+a)^2}{w_2-a}\\[1.2ex]
		-\dfrac{(w_2-a)(w_2-a-1)^2}{w_2-a-(w_1+a)^2}
	\end{pmatrix}
	-
	\begin{pmatrix}
		a\\-a-1
	\end{pmatrix}.
\end{align*}
One computes that for every $n$,
\begin{align*}
	\hat{G}^n
	=\varphi^{-1}G^n\varphi.
\end{align*}	
Since  $\{G^n\}$ is equicontinuous on the compact set $K$, and the Fubini-Study metric and the Euclidean metric are equivalent, the induced family $\{\hat{G}^n\}$ is equicontinuous on $K_0$.

To obtain a contradiction, we show that $\{\hat{G}^n\}$ cannot be equicontinuous on $K_0$. Since $\hat{G}(\mathbf{0})= \mathbf{0}$, it is enough to find $\varepsilon_0>0$ such that for every $\delta>0$, there exist $z_\delta\in K_0$ with $\|z_\delta\|<\delta$ and an integer $N_\delta>0$ satisfying $\|\hat{G}^{N_\delta}(z_\delta)\|\ge\varepsilon_0$.

Note that $\hat{G}$ is holomorphic in a neighborhood of $\mathbf{0}$, and its Jacobian at $\mathbf{0}$ is
\[
M=
\begin{pmatrix}
	2&1\\2&3
\end{pmatrix}.
\]
By Taylor expansion and the Cauchy estimates (see \cite[Theorem 1.54 and Corollary 1.101]{Sc}), there exist constants $\lambda>0$ and $C>0$ such that $B(\mathbf{0},\lambda)\subset K_0$, and for every $z\in B(\mathbf{0},\lambda)$,
\[
\hat{G}(z)=Mz+R(z),\qquad \|R(z)\|\le C\|z\|^2.
\]
The eigenvalues of $M$ are $4$ and $1$.

Set $\varepsilon_0=\min\{\lambda, (2C)^{-1}\}$. For any $\delta>0$, choose $z_\delta$ to be an eigenvector of $M$ corresponding to the eigenvalue $4$ with $\|z_\delta\|=\mu$, where $\mu=\frac12\min\{\delta,\varepsilon_0\}$. Then $\|z_\delta\|<\delta$ and $z_\delta\in B(0,\varepsilon_0)\subset K_0$. Now two cases are possible.

\textbf{Case (i)} There exists $N_\delta>0$ such that $\|\hat{G}^{N_\delta}(z_\delta)\|\ge\varepsilon_0$. Then we are done.

\textbf{Case (ii)} For all $n>0$, $\|\hat{G}^n(z_\delta)\|<\varepsilon_0$. We show this leads to a contradiction.

Let $P=(M-I)/3$ and $Q=(4I-M)/3$. One checks that $P+Q=I$, $PQ=QP=O$, $P^2=P$, $Q^2=Q$, and the operator norms $\|P\|=\|Q\|=\frac{\sqrt{10}}{3}$. Define the cone $\mathcal{N}=\{z\in\mathbb{C}^2:\|Qz\|<\|Pz\|\}$.

We claim that if $z\in\mathcal{N}\cap B(0,\varepsilon_0)$, then $\hat{G}(z)\in\mathcal{N}$. Indeed, write $z=Pz+Qz$ with $\|Qz\|<\|Pz\|$. Then
\[
\|R(z)\|\le C\|z\|^2\le C\|z\|(\|Pz\|+\|Qz\|)<2C\varepsilon_0\|Pz\|\le \|Pz\|.
\]
Moreover,
\[
\hat{G}(z)=Mz+R(z)=4Pz+Qz+R(z),
\]
and hence 
\[
P\hat{G}(z)=4Pz+PR(z),\qquad Q\hat{G}(z)=Qz+QR(z).
\]
Consequently,
\begin{align}\label{PGz}
	\|P\hat{G}(z)\|\ge 4\|Pz\|-\|P\|\|R(z)\|
	> \left(4-\frac{\sqrt{10}}{3}\right)\|Pz\|,
\end{align}
and
\[
\|Q\hat{G}(z)\|\le \|Qz\|+\|Q\|\|R(z)\|
< \left(1+\frac{\sqrt{10}}{3}\right)\|Pz\|.
\]
Since $4-\frac{\sqrt{10}}{3}>1+\frac{\sqrt{10}}{3}$, we obtain $\|Q\hat{G}(z)\|<\|P\hat{G}(z)\|$, i.e. $\hat{G}(z)\in\mathcal{N}$. This proves the claim.

Now, for our chosen $z_\delta$, we have $Pz_\delta=z_\delta$ and $Qz_\delta=0$, hence $z_\delta\in\mathcal{N}\cap B(0,\varepsilon_0)$. By the claim, $\hat{G}(z_\delta)\in\mathcal{N}$. Combining with the assumption of Case (ii), we get $\hat{G}(z_\delta)\in\mathcal{N}\cap B(0,\varepsilon_0)$. An induction then yields $\hat{G}^n(z_\delta)\in\mathcal{N}\cap B(0,\varepsilon_0)$ for all $n\ge0$.

For any $z\in\mathcal{N}\cap B(0,\varepsilon_0)$, we have
\[
\|\hat{G}(z)\|\ge 4\|Pz\|-\|Qz\|-\|R(z)\|>2\|Pz\|.
\]
Applying this to $z=\hat{G}^{n-1}(z_\delta)$ gives
\[
\|\hat{G}^n(z_\delta)\|=\|\hat{G}(\hat{G}^{n-1}(z_\delta))\|>2\|P\hat{G}^{n-1}(z_\delta)\|.
\]
Iterating the estimate \eqref{PGz} yields, for $n\ge2$,
\[
\|\hat{G}^n(z_\delta)\|
>2\|P\hat{G}^{n-1}(z_\delta)\|
>2\left(4-\frac{\sqrt{10}}{3}\right)^{n-1}\|Pz_\delta\|
>2^n\mu,
\]
as $\|Pz_\delta\|=\|z_\delta\|=\mu$. Thus $2^n\mu<\|\hat{G}^n(z_\delta)\|<\varepsilon_0$ for every $n\ge2$, which is impossible for sufficiently large $n$. Therefore Case (ii) cannot hold, giving the desired contradiction. Hence $z_a\in\mathcal{J}$.\\

\textbf{Case 2.} If $z$ is a point in $L$ with at least one zero coordinate, then $z\in\mathcal{J}$.

For $z=[z_0:z_1:z_2]\in L$, if $z_1$ or $z_2$ equals zero, it follows from Proposition \ref{indeterminacy_set} that $z\in E(F)$. Hence $z\in\mathcal{J}$. 

If $z\in L$ with $z_0=0$, then $z=[0:1:-1]$. We show that $z\in\mathcal{J}$. Suppose, for contradiction, that $z$ is a Fatou point. Then, there exists a neighborhood $U\subset \mathbb{P}^2\setminus E(F)$ of $z$ such that $\{F^n\}$ is a normal family on $U$. Let $K$ be a compact subset of $U$ with $z$ in its interior. Then there exists a subsequence $\{F^{n_k}\}$ which converges uniformly on $K$, to a continuous function $f$.

For any $n$, set $z_n=\left[\frac{1}{2^n}:1+\frac{1}{2^n}:\frac{1}{2^n}-1\right]=[1:2^n+1:1-2^n]$. Since $z$ lies in the interior of $K$, and $z_n$ converges to $z$, there exists an integer $N$ such that $z_n\in K$ for all $n>N$. For each fixed $z_n$, since $\left\lvert 2^n+1\right\rvert-\left\lvert 1-2^n \right\rvert>1$, the proof of Proposition \ref{R_1_Fatou} implies that $F^{2k-1}(z_n)$ converges to $[0:0:1]$ and $F^{2k}(z_n)$ converges to $[0:1:0]$ as $k\to\infty$. Note that for $n>N$, $F^{n_k}(z_n)$ is convergent, 
and hence $F^{n_k}(z_n)\to [0:0:1]$ or $F^{n_k}(z_n)\to [0:1:0]$, as $k\to\infty$ .
This yields $f(z_n)=[0:0:1]$ or $f(z_n)=[0:1:0]$. On the other hand, since $z_n$ converges to $z$, we have for each $n_k$, $F^{n_k}(z_n)\to F^{n_k}(z)$ as $n\to \infty$. One checks that $z$ is a fixed point of $F$, and thus $F^{n_k}(z_n)\to z$ as $n\to \infty$. Letting $k\to\infty$, since $\{F^{n_k}\}$ converges uniformly to $f$, we get that $f(z_n)\to z$, which is a contradiction. This indicates that $z\in\mathcal{J}$.\\

Combining the proofs of Case 1 and Case 2, we conclude that $L\subseteq \mathcal{J}$. Hence Proposition \ref{L_in_Julia} is proved.
\zb\\

Theorem \ref{Julia_spectrum} proved that the projective spectrum $p(A_{\rho})$ and the Julia set $\mathcal{J}$ have a substantial overlap. Nevertheless, whether the two sets are actually identical is not yet known. The following corollary establishes that they do coincide on points with at least one zero coordinate.

\begin{corr}\label{commonJandP}
	For $z=[z_0:z_1:z_2]\in \mathbb{P}^2$, the following hold:\\
	\emph{(i)} $p(A_{\rho})\cap \{z_0=0\}=\mathcal{J}\cap \{z_0=0\}=\{z\in \mathbb{P}^2 : z_0=0,\ |z_1|=|z_2|\}$,\\
	\emph{(ii)} $p(A_{\rho})\cap \{z_1=0\}=\mathcal{J}\cap \{z_1=0\}=\{z\in \mathbb{P}^2 : z_1=0,\ |z_0|=|z_2|\}$,\\
	\emph{(iii)} $p(A_{\rho})\cap \{z_2=0\}=\mathcal{J}\cap \{z_2=0\}=\{z\in \mathbb{P}^2 : z_2=0,\ |z_0|=|z_1|\}$.
\end{corr}	

\pf
In view of Theorem \ref{Fatou_resolvent} (iii), we know that $\{z\in \mathbb{P}^2 : z_0=0,\ |z_1|\neq|z_2|\}$ is contained in $p^c(A_{\rho})\cap\mathcal{F}$. Hence to prove (i), we only need to show that $\{z\in \mathbb{P}^2 : z_0=0,\ |z_1|=|z_2|\}$ is contained in $p(A_{\rho})\cap\mathcal{J}$. According to Theorem \ref{Julia_spectrum}, it suffices to show that $\{z\in \mathbb{P}^2 : z_0=0,\ |z_1|=|z_2|\}\subseteq \overline{\bigcup^{\infty}_{n=0} F^{-n}(L)}$.

Any point in $\{z\in \mathbb{P}^2 : z_0=0,\ |z_1|=|z_2|\}$ can be written as $z_t=[0:1:t]$ with $|t|=1$. Since $F(z_t)=[0:t^2:-1]=[0:1:-t^{-2}]$, by induction, we have
\begin{equation}\label{converge_0_1_a}
	F^n(z_t)=
		\left[0:1:-t^{(-2)^n}\right],\qquad n\geq 1.	
\end{equation}

For any $n\geq 1$, let $t$ be a complex number such that $t^{2^n}=1$.
It follows by \eqref{converge_0_1_a} that $F^{n}(z_t)=[0:1:-1]\in L$. Hence 
\begin{equation*}
	\left\{z_t: t^{2^n}=1\right\}\subseteq F^{-n}(L).
\end{equation*}
Consequently, 
\begin{equation*}
	\bigcup_{n=1}^{\infty}\left\{z_t: t^{2^n}=1\right\}\subseteq \bigcup^{\infty}_{n=0} F^{-n}(L).
\end{equation*}
We take closure on both sides, then
\begin{equation*}
	\left\{z_t: \left\lvert t\right\rvert=1\right\}
	=\overline{\bigcup_{n=1}^{\infty}\left\{z_t: t^{2^n}=1\right\}}
	\subseteq \overline{\bigcup^{\infty}_{n=0} F^{-n}(L)}.
\end{equation*}
This proves (i).  

For (ii), let $z=[z_0:0:z_2]\in \mathbb{P}^2$ with $\left\lvert z_0\right\rvert=\left\lvert z_2\right\rvert$. Proposition \ref{indeterminacy_set} implies that $z\in E(F)$, and thus $z\in \mathcal{J}$. Note that Corollary \ref{corr_E_in_p} has showed that $E(F)\subseteq p(A_{\rho})$. Hence $z\in p(A_{\rho})\cap\mathcal{J}$. 
If $z=[z_0:0:z_2]\in \mathbb{P}^2$ with $\left\lvert z_0\right\rvert\neq\left\lvert z_2\right\rvert$, one checks that $z$ belongs to the set $R_0\cup R_2$ in Theorem \ref{Fatou_resolvent}. Then it follows that $z\in p^c(A_{\rho})\cap\mathcal{F}$. Combining the two inclusions yields (ii).

The proof of (iii) is similar to that of (ii). Thus the corollary follows.\zb\\

Taking complements in Corollary \ref{commonJandP} immediately yields that the Fatou set and the projective resolvent set also coincide on points with at least one zero coordinate. This is the following corollary.

\begin{corr}\label{commonFandPc}
	For $z=[z_0:z_1:z_2]\in \mathbb{P}^2$, the following hold:\\
	\emph{(i)} $p^c(A_{\rho})\cap \{z_0=0\}=\mathcal{F}\cap \{z_0=0\}=\{z\in \mathbb{P}^2 : z_0=0,\ |z_1|\neq |z_2|\}$,\\
	\emph{(ii)} $p^c(A_{\rho})\cap \{z_1=0\}=\mathcal{F}\cap \{z_1=0\}=\{z\in \mathbb{P}^2 : z_1=0,\ |z_0|\neq |z_2|\}$,\\
	\emph{(iii)} $p^c(A_{\rho})\cap \{z_2=0\}=\mathcal{F}\cap \{z_2=0\}=\{z\in \mathbb{P}^2 : z_2=0,\ |z_0|\neq |z_1|\}$.\\
\end{corr}

\section{Conclusions and further remarks}

For holomorphic dynamical systems in several variables \cite{Fo, MNTU}, the complexity of the Fatou and Julia sets makes it uncommon for a system to admit a full description of both. In this paper, since the dynamical system \eqref{eqmain} comes from a self-similar group, namely the Basilica group, we establish a novel correspondence between the Fatou set and the projective resolvent set, and between the Julia set and the projective spectrum. Although a complete identification is still lacking, we demonstrate that these sets agree on a large common subset. This naturally leads to the following questions:

\begin{que}
	Are the three sets $\mathcal{F}$, $p^c(\mathcal{B}_{\rho})$ and $\bigcup^{\infty}_{n=0} F^{-n}\left(\bigcup_{k=0}^2R_k\right)$ pairwise equal? If not, are any two of them equal?
\end{que}

\begin{que}
	Are the three sets $\mathcal{J}$, $p(A_{\rho})$ and $\overline{\bigcup^{\infty}_{n=0} F^{-n}(L)}$ pairwise equal? If not, are any two of them equal?
\end{que}

Numerical simulations further suggest that, in both cases, the common subset identified above is in fact the entire intersection, thereby indicating an affirmative answer to both questions. Note that, in the proofs of the theorems, it is relatively easy to verify the inclusion of a given subset in the projective spectrum (resp. projective resolvent set), but considerably harder to establish the corresponding inclusion in the Julia set (resp. Fatou set). This suggests that the projective spectrum can serve as a useful indicator of the corresponding Julia set. Our subsequent goal is therefore to obtain a complete correspondence. We defer further details to a separate paper.\\

\noindent\textbf{Acknowledgments.} The work in this paper is partially supported by the National Natural Science Foundation of China Grant No. 12271090, Jiangsu Provincial Scientific Research Center of Applied Mathematics Grant No. BK20233002.\\

\vspace{5mm}
\end{document}